\documentclass[12pt,twoside]{amsart}

\usepackage{latexsym,amsmath,amscd,array}
\input amssym.def
\input amssym.tex

\swapnumbers

\parskip=\medskipamount
\newtheorem{theorem}{Theorem}[section]
\newtheorem{lemma}[theorem]{Lemma}
\newtheorem{definition}[theorem]{Definition}
\newtheorem{proposition}[theorem]{Proposition}
\newtheorem{corollary}[theorem]{Corollary}
\newtheorem{comment}[theorem]{Comment}
\newtheorem{remark}[theorem]{Remark}

\numberwithin{equation}{section}

\def\p{\m {\it Proof. }}
\def\m{\medskip}

\def\C{\mathbb C}
\def\R{\mathbb R} 
\def\L{\mathbb L}
    
\def\eps{\varepsilon}

\def\Im{\operatorname {Im}}

\def\dim{\operatorname {dim}}
\def\hr{_{\operatorname {hr}}}
\def\rk{\operatorname {rank}}
\def\span{\operatorname {span}}

\def\osi{\Omega_{\operatorname {sympl}}}

\def\ga{\alpha}
\def\gb{\beta}

\def\gd{\delta}
\def\eps{\varepsilon}

\def\gl{\lambda}

\def\gL{\Lambda}

\def\asl{\mathfrak s\mathfrak l\,(2,\C)}

\newcommand{\ba}{\begin{array}}
\newcommand{\ea}{\end{array}}
\newcommand{\hb}[1]{\hbox{#1}}

\newcommand{\op}{\oplus}

\newcommand{\SSS}{{\mathcal S}}

\newcommand{\g}{{\mathfrak g}}

\begin{document}
\title{Symplectically Harmonic Cohomology of 
Nilmanifolds}
\author{R. IB\'A\~NEZ}
\author{YU. RUDYAK}
\author{A. TRALLE}
\author{L. UGARTE}
\address{R. Ib\'a\~nez, Departamento de Matem\'aticas, Facultad de
Ciencias, Universidad del Pais Vasco, Apdo. 644, 48080 Bilbao, Spain,}
\email {mtpibtor@lg.ehu.es}
\address {Yu. Rudyak, Department of Mathematics, University of Florida, 358 
Little Hall, PO Box 118105, Gainesville, FL 32611-8105, USA}
\email{rudyak@math.ufl.edu} 
\email{rudyak@mathi.uni-heidelberg.de}
\address{A. Tralle,  Department of Computer Science, The College of Economics  
and Computer Science TWP, Wyzwolenia str. 30, 10-106 Olsztyn, Poland}
\email{tralle@matman.uwm.edu.pl}
\address{L. Ugarte, Departamento de Matem\'aticas, Facultad de
Ciencias. Universidad de Zaragoza, 50009 Zaragoza, Spain}
\email{ugarte@posta.unizar.es}

\maketitle

\section{Introduction}

This paper can be considered as an extension to our paper \cite{IRTU}. Also, it 
contains a brief survey of recent results on symplectically harmonic cohomology. 
Below $H^k(M)$ always means the de Rham cohomology of a smooth manifold $M$.

Let $(M^{2m}, \omega)$ 
be symplectic manifold, and let $\Omega^k(M)$ be the 
space 
of all $k$-forms on $M$. In \cite{L,Br} it was introduced an operator 
$$
*: \Omega^k(M) \to \Omega^{2m-k}(M), \quad **=1
$$
which is a symplectic analog of the well-known de Rham--Hodge $*$-operator on 
oriented Riemannian manifolds: one should use 
the symplectic form instead of the Riemannian metric. Going further, one can 
define operator $\gd=\pm *d*, \gd^2=0$. The form $\ga$ is called {\it 
symplectically harmonic} if $d\ga=0=\gd\ga$. However, unlike de Rham--Hodge 
case, there exist simplectically harmonic forms which are exact. Because of 
this, Brylinski \cite{Br} defined the symplectically harmonic cohomology 
$H^*\hr(-)$ by setting
$$
H^k\hr(M)=H^k\hr(M,\omega):=\Omega^k\hr(M)/(\Im\, d\cap
\Omega^k_{\hr}(M))
$$
where $\Omega^k\hr(M)$ is the space of all simplectically harmonic $k$-forms. We 
set $h_k=h_k(M,\omega) =\dim H^k\hr(M,\omega$. Since $H^k\hr(M) \subset H^k(M)$, 
we conclude that $h_k\le b_k$. Brylinski proved that $h_k=b_k$ if $M$ is a 
K\"ahler manifold. However, this equality does not hold in general: Mathieu 
\cite{Ma} proved that the equalities $h_k=b_k, k=0,1, \ldots, m$, hold if and 
only if $M$ has the Hard Lefschetz Property. Since there are many simplectic 
manifolds without K\"ahler structure (see e.g. a survey \cite{TO}), we have many 
manifolds with $h_k < b_k$.

\m The next step is to ask whether $h_k(M, \omega)$ can vary with respect to 
$\omega$. According to Yan \cite{Y}, the following question was posed by Boris 
Khesin and Dusa
McDuff. 

\m {\bf Question:} Do there exist compact manifolds $M$ which possess a 
continuous
family $\omega_t$ of symplectic forms such that $h_k(M,
\omega_t)$ varies with respect to $t$?

This question, according to Khesin, is probably related to group
theoretical hydrodynamics and geometry of diffeomorphism groups. Some
indirect indications for an existence of such relations can be found
in \cite{AK}.

\m In \cite{IRTU} we named such manifolds {\it flexible} and proved that there 
are several flexible manifolds among 6-dimensional nilmanifolds. Notice that Yan 
\cite{Y} constructed a 4-dimensional flexible manifold (in fact, his arguments 
need a small correction, see \cite{IRTU}.) Sakane and Yamada \cite{SY} also 
found some flexible manifolds among nilmanifolds.

\m One can prove that, for every 6-dimensional manifold, $h_i=b_i$ for 
$i=0,1,2,6$. In \cite{IRTU} we computed $h_4$ and $h_5$ for 6-dimensional 
nilmanifolds, and in this paper, using a result of Yamada \cite{Ym}, we compute 
$h_3$. As a corollary, we prove that there are exactly 10 flexible manifolds 
among 6-dimensional nilmanifolds.

\m It is interesting to mention the following phenomenon. For generic symplectic 
form $\omega$ on $6$-dimensional manifolds, each of the numbers $h_i(M,\omega), 
i =4,5$ takes the maximal value. In other words, the set of symplectic forms 
$\omega$ with the maximal value of $h_i(\omega), i=4,5$ is open and dense in the 
set of all symplectic forms, see Corollary \ref{strat}. On the other hand, if 
$h_5$ does not depend on $\omega$ then, for generic symplectic form $\omega$ on 
$6$-dimensional manifolds, the number $h_3(\omega)$ takes {\it minimal} value, 
see Corollary \ref{strat2}.

Finally, we notice that all known examples of flexible manifolds are non-simply 
connected. So, it would be interesting to have  examples of simply connected 
flexible 
manifolds.

\section{Some Operators on Symplectic Vector Spaces}\label{local}

Let $V^{2m}$ be a real vector space of dimension $2m$, let $V^*$ be the 
adjoint space, and let $\gL^k(V^*)$ be the $k$-th exterior power of $V^*$. In 
other words, $\gL^k(V^*)$ is the space of skew-symmetric $k$-forms on $V$. A 
{\it symplectic form} on $V$ is a non-degenerate 2-form $\eta$, that is, $\eta 
\in \gL^2(V^*)$ such that $\eta^m\ne 0$. A {\it symplectic vector space} is a 
vector space with a given symplectic form.  

The form $\eta$ yields a linear map $V \to V^*$ of the form $v \mapsto 
i(v)\eta$. Here $i$ is the contraction (interior product). This map is an 
isomorphism since $\eta$ is non-degenerate. More 
generally, we have an isomorphism
$$
\CD
\mu=\mu_k: \{k\text{-vectors}\} \to \{k\text{-forms}\},
\endCD
$$
where $\mu_k$ extends $\mu_1$ such that it respects the exterior 
multiplication.
Let $\Pi$ be the bivector dual to the 
form $\eta$. We define the operator

\begin{equation}\label{star-br}
*:\gL^k(V^*)  \longrightarrow \gL^{2m-k}(V^*)
\end{equation}

by requiring

\begin{equation}\label{*-poisson}
\gb \wedge (*\ga)=\gL^k(\Pi)(\gb, \ga)\frac{\eta^m}{m!}
\end{equation}

for every two $k$-forms $\ga$ and $\gb$.

\m
The operator $*$ has the following properties. Let $(V_1, \eta_1)$ and $(V_2, 
\eta_2)$ be two symplectic vector spaces, and let $*_1$ and $*_2$ be the 
related 
$*$-operators. Let $p_i: V_1\times V_2 \to V_i, i=1,2$ be the projections. We 
equip $V_1 \times V_2$ with the symplectic form $p_1^*\eta_1+p_2^*\eta_2$ and 
denote the related $*$-operator just by $*$. Let $\ga_1$ be a $k_1$-form on 
$V_1$, and let $\ga_2$ be a $k_2$-form on $V_2$. We set
$$
\ga_1\boxtimes \ga_2= (p_1^*\ga_1)\wedge(p_2^*\ga_2).
$$

\begin{proposition}[\cite{Br}]
\label{prod}
$$
*(\ga_1\boxtimes \ga_2)=(-1)^{k_1k_2}(*_1\ga_1)\boxtimes(*_2\ga_2).
$$
\qed
\end{proposition}

\begin{proposition}[\cite{Br}]
\label{m=1}
If $\dim V =2$ then $*\ga=-\ga$ for every $1$-form $\ga$. Furthermore, 
$*\eta=1$ 
and $*(1)=\eta$.
\qed
\end{proposition}

In fact, these properties allow us to evaluate the $*$-operator on any 
symplectic vector space.
 In particular, one can prove that $**\ga=\ga$.

\m
We define the maps 
$$
L=L_{\eta}: \gL^k (V^*) \to \gL^{k+2}(V^*), \quad L(\ga)=\ga \wedge \eta
$$
and
$$
L^*=-*L*: \gL^k (V^*) \to \gL^{k-2}(V^*).
$$
Yan \cite{Y} proved that $L^*=i(\Pi)$. Finally, we define the map 
$$
A: \gL^*(V^*) \to   \gL^*(V^*), \quad A = \sum (m-k) \pi_k  
$$
where $\pi_k: \gL^*(V^*)  \to \gL^k(V^*) $ is the natural projection.

\m Certainly, we can regard $L$ and $L^*$ as maps $\gL^*(V^*) \to   
\gL^*(V^*)$. 
It is easy to see that the following relations hold, see e.g. \cite{Y}:

\begin{equation}\label{sympl-sl}
[L^*,L]=A, \quad [A,L]=-2L, \quad [A,L^*]=2L^*.
\end{equation} 

\begin{remark}\rm
The symplectic $*$-operator was originally considered by Libermann \cite{L}, 
see also \cite{LM}, as
\begin{equation}\label{lib}
*:\gL^k(V^*)  \longrightarrow \gL^{2m-k}(V^*), \quad 
*\ga=i(\mu_k^{-1}(\ga))\frac{\eta^m}{m!}.
\end{equation}

She also introduced the operator $L^*=-*L*$. In order to see that the 
$*$-operators from \eqref{*-poisson} and \eqref{lib} coincide, it suffices to 
prove that the Libermann's $*$-operator has the properties from Propositions 
\ref{prod} and \ref{m=1}. This can be done immediately.
\end{remark}

\m Consider the Lie algebra $\asl$. It is generated over $\C$ by the matrices
$$
X=\left(\begin{array}{cc}
0&1\\
0&0
\end{array}\right)
\quad
Y=\left(\begin{array}{cc}
0&0\\
1&0
\end{array}\right)
H=\left(\begin{array}{cc}
1&0\\
0&-1
\end{array}\right).                                       
$$

It is easy to see that

\begin{equation}\label{sl}
[X,Y]=H,\quad [H,X]=2X, \quad [H,Y]=-2Y.
\end{equation}

The proof of the folowing well-known lemma can be found e.g. in \cite{GH, W}.

\begin{lemma}
\label{sl-rep}
Let $V$ be a finite dimensional complex vector space which is a space of a 
representation of $\asl$. Then all eigenvalues of $H$ are integers. Let $V_k$ be 
the eigenspace of $H$ with respect to 
eigenvalue $k$. Then
\begin{equation*}
Y^k: V_{-k} \rightarrow V_k \text{ and } X^k: V_k \rightarrow V_{-k}
\end{equation*} 
are isomorphisms.
\qed 
\end{lemma}

Since 
relations \eqref{sympl-sl} have the same form as relations \eqref{sl},  we can 
equip $\gL^*({V_{\C}}^*)$, where $V_{\C}$ is the complex vector
space $V\otimes \C$, with a structure of 
$\asl$-module via the representation
$$
X \mapsto L^*, \quad Y \mapsto L, \quad H \mapsto A.
$$ 

Clearly, in this representation $V_k=\gL^k({V_{\C}}^*)$. Now, since $L$ is
a real operator, Lemma \ref{sl-rep} implies that

\begin{equation}\label{finite}
L^k: \gL^{m-k}(V^*) \rightarrow \gL^{m+k}(V^*)
\end{equation}

is an isomorphism.

\section{Symplectically Harmonic Forms}

Given a smooth manifold $M$, a {\it symplectic form} on $M$ is a closed 2-form 
$\omega$ on $M$ such  that, for every $p\in M$, the form $\omega|_{T_pM}$ is a 
symplectic form on $T_pM$. A symplectic manifold is a manifold with a fixed 
symplectic form. We denote by $\Pi$ the bivector field dual to $\omega$. 

\m Let $\Omega^k(M)$ be the space of all $k$-forms on $M^{2m}$. Given a 
symplectic form $\omega$ on $M$, we can introduce the $*$-operator
$$
\CD
*: \Omega^k(M) @>>> \Omega^{2m-k}(M)
\endCD
$$
just as we did it in \eqref{*-poisson}. Furthermore, we have the operator
$$
L=L_{\omega}: \Omega^k (M) \to \Omega^{k+2}(M), \quad L(\ga)=\ga \wedge \omega
$$
and
$$
L^*=-*L*: \Omega^k (M)\to \Omega^{k-2}(M).
$$
Finally, there is an operator
$$
A: \Omega^*(M) \to   \Omega^*(M), \quad A = \sum (m-k) \pi_k  
$$
and the relations \eqref{sympl-sl} hold. However, in order to get an analog of 
\eqref{finite}, we need to have a generalization of Lemma \ref{sl-rep} for 
infinite dimensional vector spaces. Following Yan \cite{Y}, we say that an 
$\asl$-representation $V$ (not necessary finitely dimensional) is of {\it finite 
$H$-spectrum} if

\begin{enumerate}
\item[(1)] $V$ can be decomposed as the direct sum of eigenspaces of $H$; 
\item[(2)] $H$ has only finitely many distinct eigenvalues.
\end{enumerate}

Yan proved that for $\asl$-modules of finite $H$-spectrum all the eigenvalues 
of $H$ are integers. Furthermore, an analog of Lemma \ref{sl-rep} holds for 
$\asl$-modules of finite $H$-spectrum, i.e., we have the following fact.

\begin{lemma}
\label{finite-h}
Let $V$ be an $\asl$-module of finite $H$-spectrum. Let $V_k$ be the eigenspace 
of $H$ with respect to eigenvalue $k$. Then
\begin{equation*}
Y^k: V_{-k} \rightarrow V_k \text{ and } X^k: V_k \rightarrow V_{-k}
\end{equation*} 
are isomorphisms.
\qed 
\end{lemma}

It is easy to see that $\Omega^*(M,\C)$ turns out to be  an 
$\asl$-module of finite $H$-spectrum. Now, asserting as in Section
\ref{local}, we conclude that
$$
\CD
L^k: \Omega^{m-k}(M) @>>> \Omega^{m+k}(M)
\endCD
$$
is an isomorphism.

For every $k$ we introduce the operator 
$$
\CD
\gd: \Omega^k(M) @>>>\Omega^{k-1}(M), \quad \gd(\ga)=(-1)^{k+1}*d(*\ga).
\endCD
$$
It turns out to be that $\gd=[i(\Pi),d]$, 
see \cite{Br}.

\begin{remark}\rm
The operator $\gd=-*d*$ was also considered by Libermann (see \cite{LM}). 
Koszul \cite{K} introduced the operator $\gd=[d,i(\Pi)]$ for Poisson manifolds. 
Brylinski \cite{Br} proved that these operators coincide.
\end{remark}

\begin{definition}[\cite{Br}]
\rm
A form $\ga$ on a symplectic manifold $(M, \omega)$ is called {\it 
symplectically harmonic} 
if $d\ga =0=\gd \ga$.
\end{definition}

We denote by $\Omega^k\hr(M)$ the linear space of symplectically harmonic 
$k$-forms. It is clear that $\Omega^k\hr(M,\C)$ is an $\asl$-submodule of 
$\Omega^k(M,\C)$. 
Thus, Lemma \ref{finite-h} yields the following result.

\begin{theorem} 
[\cite{Y}]
\label{iso}
The map
$$
\CD
L^k: \Omega^{m-k}\hr(M) @>>> \Omega^{m+k}\hr(M)
\endCD
$$
is an isomorphism.
\qed
\end{theorem}

\m Unlike the Hodge theory, there are non-zero exact symplectically harmonic 
forms. 
Now, following Brylinski \cite{Br}, we define symplectically harmonic 
cohomology 
$H^*\hr(M,\omega)$ by setting
$$
H^k\hr(M,\omega)=\Omega^k\hr(M)/(\Im d \cap \Omega^k\hr(M))
$$
and 
$$
h_k=h_k(M,\omega)=  \dim  H^k\hr(M,\omega).
$$
So, $H^k\hr(M,\omega) \subset H^k(M)$.

\m
Notice that Theorem \ref{iso} implies the following corollary.
\begin{corollary}\label{epi}
$$
\CD
L^k: H^k\hr(M) @>>> H^{2m-k}\hr(M)
\endCD
$$
is an epimorphism. In particular, $h_{m-k}(M) \ge h_{m+k}(M)$.
\qed
\end{corollary}

We set 
$$
P^{m-k}(M,\omega)=\{a\in H^{m-k}(M)\bigm| L^{k+1}a=0\}.
$$
Yan \cite{Y} proved that $P^{m-k}(M,\omega) \subset H^{m-k}\hr(M,\omega)$.
The following result allows us to describe the groups $H^k\hr(M,\omega)$. 

\begin{theorem}\label{harm}
For every $k\ge 0$ we have
\begin{eqnarray}
H^{m-k}\hr(M)&=&P^{m-k}(M)+L(H^{m-k-2}\hr(M)) \subset 
H^{m-k}(M);\nonumber\\
H^{m+k}\hr(M)&=&\Im\{L^k:H^{m-k}\hr(M) \longrightarrow H^{m+k}(M)\} \subset 
H^{m+k}(M).\nonumber
\end{eqnarray}
\qed
\end{theorem}

The first equality is proved in \cite[Lemma 4.3]{Ym}, the second equality is 
proved in 
\cite[Corollary 2.4]{IRTU}.

\begin{corollary}
[\cite{Ym}]
\label{cohom}
If symplectic forms $\omega$ and $\omega'$ belong to the same cohomology 
class, then $H^*\hr(M, \omega)=H^*\hr(M, \omega')$.
\qed
\end{corollary}

We need also the following fact. It can be deduced from Theorem \ref{harm} 
directly, or see \cite{Y, IRTU}.

\begin{proposition}
\label{small}
$h_k=b_k$ for $k=0,1,2$.
\qed
\end{proposition}

\section{Flexibility}

\begin{definition}\rm
We say that a smooth (closed) manifold $M$ is {\it $k$-flexible} if $M$ 
possesses a continuous family of symplectic forms $\omega_t, t\in [0,1]$ such 
that $h_k(M, \omega_0) \ne h_k(M,\omega_1)$. We also say that a manifold is 
flexible if it is $k$-flexible for some $k$. 
\end{definition}

\m 
Certainly, the existence of two symplectic forms $\omega_1, \omega_2$ with 
$h_k(\omega_1)\ne h_k(\omega_2)$ is necessary for $k$-flexibility, but is not 
sufficient in general. However, as we will see below, this condition is 
sufficient for $k$-flexibility if $k\ge \dim M-2$. It was proved in \cite{IRTU}, 
here we arrange this proof in more explicit way.

\begin{lemma}\label{linalg}
Let $\L$ be the space of all linear maps $\mathbb R^m \to \mathbb R^m$. Fix any 
linear map $D:\mathbb R^m \to \mathbb R^l$ and a positive integer $k$. Let 
$A,B\in \L$ be two linear maps such that $\rk\, DA^k <\rk\, 
DB^k$. 
Then the set
$$
\Lambda=\{\lambda \in \R \bigm | \rk\,D(A+\lambda B)^k >\rk\,DA^k\}
$$
is an open and dense subset of $\R$.
\end{lemma}

\p The set $\R \setminus \Lambda$ is an algebraic
 subset of $\R$, because the rank of a matrix is equal to the
 order of the largest non-zero minor. So, it suffices to prove that $\Lambda\neq
 \emptyset$. But, clearly, 
$$
\rk\, D(B+\mu A)^k\ge \rk\, DB^k > \rk\, DA^k 
$$ 
for $\mu$ small enough, and
so $\Lambda\neq\emptyset$.
\qed

\begin{corollary}\label{form}
Let $\ga, \gb$ be two closed $2$-forms on a manifold $M$. Suppose that
$$
\rk\,L^k_{\ga}< \rk\,L^k_{\gb}
$$
where
$$
L^k_{\ga}, L^k_{\gb}: H^{m-k}(M) \to H^{m+k}(M).
$$
Then the set 
$$
\Lambda=\{\lambda \in \R \bigm | \rk(L^k_{\ga}+\lambda L^k_{\gb})> 
\rk\,L^k_{\ga}\}
$$
is an open and dense subset of $\R$. Furthermore, for every $b\le \rk L_{\gb}$ 
the set of closed $2$-forms
$$
U=\{u\mid \rk\,\{L^k_{u}: H^{m-k}(M) \to H^{m+k}(M)\} \ge b\}
$$
is an open and dense subset of the space of closed $2$-forms.
\end{corollary} 

\p Notice that the map
$$
L^k_{\ga}: H^{m-k}(M) \to H^{m+k}(M)
$$
can also be written as
$$
\CD
H^*(M) @>L^k>> H^*(M) @>P>>H^{m+k}(M)
\endCD
$$
where $P$ is the obvious projection. Now, applying Lemma \ref{linalg} with 
$A=L_{\ga}, B=L_{\gb}$ and $D=P$, we get the desired result about $\gL$. 
Furthermore, $U$ is open by general reasons (small perturbation does not 
decrease the rank). We prove that $U$ is dense. Take any 2-form $\gamma$ with 
$\rk\, L^k_{\gamma}<m$. Then, by what we said above, there exists arbitrary 
small $\gl$ with  $\rk\, L^k_{\gamma+\gl\gb}=\rk L^k_{\gamma}+\gl L^k_{\gb}\ge 
m$. Thus, $U$ is dense.
\qed
   
\begin{corollary}\label{flex1}
Let $\omega_0$ and $\omega_1$ be two symplectic forms on a manifold
$M^{2m}$. Suppose that, for some $k>0$, $h_{m-k}(\omega_0)=
h_{m-k}(\omega_1)$, but $h_{m+k}(\omega_0)< h_{m+k}(\omega_1)$. Then,
for every $\eps>0$, there exists $\lambda\in (0,\eps)$ such that
$h_{m+k}(\omega_0+\lambda
\omega_1)> h_{m+k}(\omega_0)$. Moreover, $M$ is 
flexible
provided that it is closed.
\end{corollary}

\p 
Because of Theorem \ref{harm},
$$
h_{m+k}=\rk \{L^k: H^{m-k}\hr(M) \to H^{m+k}(M)\}
$$
Now the existence of $\gl$ follows from Corollary \ref{form}. To prove the 
flexibility of $M$, take the above $\gl$ 
so small that $\omega_0+t\omega_1$ is a symplectic form for $t\in [0,\lambda]$ 
small enough. Now,
we set $\omega_t=\omega_0+t\omega_1, t\in [0,\lambda]$ and 
conclude that $h_{m+k}(\omega_0) < h_{m+k}(\omega_{\lambda})$.
\qed

\begin{corollary}\label{flex2}
Let $(M^{2m}, \omega_0)$ be a closed
symplectic manifold. Fix any $k$ with $0<k<m$ and suppose that
$h_{m-k}(M,\omega)= b_{m-k}(M)$ for every symplectic form $\omega$ on
$M$. Furthermore, suppose that there exists $x\in H^2(M)$ such that
$$
\rk\,\{L^k_x: H^{m-k}(M) \to H^{m+k}(M)\} > h_{m+k}(M,\omega_0).
$$
Then $M$ is flexible.
\end{corollary}

\p Take a closed 2-form $\alpha$ which represents $x$. Then
$\omega_0+t\alpha$ is a symplectic form for $t$ small enough.
Using Theorem \ref{harm} and Corollary \ref{form} and asserting as in the 
proof of \ref{flex1}, we conclude that there exists a small $\gl$ with 
$h_{m+k}(\omega_0+\gl\ga)>h_{m+k}(\omega_0)$. Now the result follows from 
Corollary
\ref{flex1}.
\qed

\begin{corollary}\label{flex3}
Let $M^{2m}$ be a closed smooth manifold. Suppose that there are two 
symplectic forms $\omega_1,  \omega_2$ on $M$ such that 
$h_{2m-k}(M,\omega_1)\ne h_{2m-k}(M,\omega_2)$ for $k=1$ or $2$. Then $M$ is 
$k$-flexible.
\end{corollary}

\p This follows directly from Corollary \ref{flex2} and Proposition 
\ref{small}.
\qed

The above results imply the following fact noticed in \cite{IRTU}. Set
$$
\osi (M)=\{\omega\in \Omega^2(M)\bigm|\omega\ \hbox{is a symplectic form on}\ 
M\}
$$
and define $\Omega(b,k)=\{\omega\in \osi\bigm|h_k(M,\omega)=b\}$.

\begin{corollary}
\label{strat}
Let $M^{2m}$ be a manifold that admits a symplectic structure. Suppose that, 
for 
some $k>0$, $h_{m-k}(M,\omega)$
does not depend on the symplectic structure $\omega$ on $M$. Then the following
three conditions are equivalent:
\par {\rm (i)} the set $\Omega (b, m+k)$ is open and dense in $\osi(M)$;
\par {\rm (ii)} the interior of the set $\Omega (b, m+k)$ in $\osi(M)$ is 
non-empty;
\par {\rm (iii)} the set $\Omega (b, m+k)$ is non-empty and
$h_{m+k}(M,\omega)\leq b$ for every $\omega\in \osi(M)$.
\end{corollary}

\p (i) $\Rightarrow$ (ii). Trivial.
\par (ii) $\Rightarrow$ (iii). Suppose that there exists $\omega_0$ with
$h_{m+k}(M, \omega_0)> b$. Take $\omega$ in the interior of
$\Omega(b,m+k)$. Then, in view of Corollary \ref{form}, there exists an
arbitrary small $\lambda$ such that $h_{m+k}(\omega+\lambda
\omega_0)>b$, i.e. $\omega$ does not belong to the interior of
$\Omega (b,m+k)$. This is a contradiction.
\par (iii) $\Rightarrow$ (i). Notice that 
$$
\Omega (b, m+k)=
\{\omega\in \osi\mid \rk\,\{L^k_{\omega}: H^{m-k}(M) \to 
H^{m+k}(M)\}\ge b\}.
$$
Now the result follows from Corollary \ref{form}.
\qed

So, the family $\{\Omega(b,m+k)|b=0,1, \ldots \}$ gives us a
stratification of $\osi (M)$ where the maximal stratum is open and
dense.

\section{Nilmanifolds}

Given a Lie algebra $\g$, we set $\g^0=\g$ and $\g^r=[\g, \g^{r-1}]$. The Lie 
algebra $\g$ is 
called {\it nilpotent} if $\g^r=0$ for some $r$. The maximal $s$ such that 
$\g^s 
\ne 0$ is 
called the {\it step length} of the nilpotent Lie algebra $\g$. 

A Lie group $G$ is called {\it nilpotent} if its Lie algebra is nilpotent. A 
{\it 
nilmanifold} is defined to be a closed manifold $M$ of the form $G/\Gamma$ 
where 
$G$ is a 
simply connected nilpotent group and $\Gamma$ is a discrete subgroup of $G$. It 
is well known 
that $\Gamma$ determines $G$ and is determined by $G$ uniquely up to 
isomorphism 
(provided 
that $\Gamma$ exists), \cite{M,R}.

Three important facts in the study of compact nilmanifolds are (see
\cite{TO}):
\begin{enumerate} \item Let $\g$ be a nilpotent Lie algebra with
structural constants $c^{ij}_k$ with respect to some basis, and let 
$\{\alpha_1,\ldots
,\alpha_n\}$ be the dual basis of $\g^{\ast}$. Then in the Chevalley--Eilenberg 
complex
$(\Lambda^{\ast}\g^{\ast},d)$ we have
\begin{equation}\label{cte}
d\alpha_k=\sum_{1\leq i<j<k} c^{ij}_k
\alpha_i\wedge \alpha_j.
\end{equation}
\item Let $\g$ be the Lie algebra of a simply connected nilpotent Lie group 
$G$. 
Then, by Malcev's theorem \cite{M}, $G$ admits a lattice if
and only if $\g$ admits a basis such that all the structural constants
are rational.
\item By Nomizu's theorem, the Chevalley--Eilenberg complex
$(\Lambda^{\ast}\g^{\ast},d)$ of $\g$ is quasi-isomorphic to the de Rham
complex of $G/\Gamma$. In particular,
\begin{equation}\label{nom}
H^{\ast}(G/\Gamma)\cong H^{\ast}(\Lambda^{\ast}\g^*,d)
\end{equation}
and any cohomology class $[a]\in H^k(G/\Gamma)$
contains a homogeneous representative $\alpha$. Here we call the form $\alpha$
homogeneous if the pullback of $\alpha$ to $G$ is left invariant.
\end{enumerate}

These results allows us to compute cohomology invariants of nilmanifolds in 
terms of the Lie algebra $\g$, and this simplifies calculations. For example, 
Yamada \cite[Theorem 3]{Ym} proved the following theorem.

\begin{theorem}
Let $M^{2m}=G/\Gamma$ be a nilmanifold of the step length $r+1$. Then for every 
symplectic structure on $M$ we have
$$
h_1(M)-h_{2m-1}(M)\ge \dim \g^r.
$$
Moreover, if $M$ has the step length $2$ then 
$$
h_1(M)-h_{2m-1}(M)=\dim[\g,\g].
$$
\qed
\end{theorem}

This theorem yields the following corollary.

\begin{corollary}
If $M^{2m}$ has the step length $2$ then 
$$
h_{2m-1}(M)=2(b_1(M)-m).
$$
\end{corollary}

\p Notice that $b_1(M)=\dim\g-\dim[\g,\g]=2m-\dim[\g,\g]$. Now,
$$
h_{2m-1}(M)=h_1(M)-\dim[\g,\g]=b_1(M)-\dim[\g,\g]=2(b_1(M)-m).
$$
\qed

\m
Yan \cite{Y} proved that there are no flexible 4-dimensional nilmanifolds. Now 
we describe all flexible 6-dimensional nilmanifolds. In \cite{IRTU} 
we computed the harmonic numbers $h_4$ and $h_5$.  Here, using Theorem 
\ref{harm}, we compute $h_3$ (see below). 

\m
The table below extends the table from \cite{IRTU}. Namely, we added the column 
which contains $h_3$. 

\m We should also mention that the table from \cite{IRTU} used the 
classification of nilpotent Lie algebras given by Salamon~
\cite{Sa}.The last one, in turn, is based on the Morozov classification of
$6$-dimensional nilpotent Lie algebras \cite{OV}.

In the table Lie algebras appear lexicographically with respect to the triple
$(b_1,b_2,6-s)$ where $s$ is the step length. The first two columns
contain the Betti numbers $b_1$ and $b_2$ (notice that
$b_3=2(b_2-b_1+1)$ because of the vanishing of the Euler
characteristic). The next column contains $6-s$.
\par The fourth column contains the description of the structure of the
Lie algebra by means of the expressions of the form (\ref{cte}) in the
Chevalley-Eilenberg complex. It means that, say,
for the compact nilmanifold $M$ from the second row, there exists a
basis $\{\alpha_i\}_{i=1}^6$ of homogeneous $1$-forms on $M$ such that
\begin{eqnarray*}
d \alpha_1=0=d \alpha_2,&&d \alpha_3=\alpha_1 \wedge \alpha_2,\quad d
\alpha_4=\alpha_1 \wedge \alpha_3,\\
 d \alpha_5=\alpha_1 \wedge
\alpha_4,&&d \alpha_6=\alpha_3 \wedge
\alpha_4 + \alpha_5 \wedge \alpha_2.
\end{eqnarray*}
\par The column headed $\oplus$ indicates the dimensions of the
irreducible subalgebras in case $\g$ is not itself irreducible.
\par The next columns show the dimensions
$h_k$ for $k=3,4,5$. So, the column, say, $h_3$ contains all possible
values of $h_3(M,\omega)$ which appear when $\omega$ runs over all
symplectic forms on $M$. The sign ``--'' at a certain row means that
the corresponding Lie algebra (as well as the compact nilmanifold)
does not admit a symplectic structure.

\par For completeness, in the last columns we list the
dimension $\dim_{\R} \SSS(\g)$ of the moduli space of symplectic structures.

\vfil\eject\begin{center}
{\large\bf Six-dimensional real nilpotent Lie algebras}\vspace{5pt}

\def\no{\hb{\bf--}}
{\extrarowheight2pt
\begin{tabular}{|c|c|c|c|c|c|c|c|c|}
\hline 
$b_1$&$b_2$&$6\!-\!s$& Structure &$\op$& $h_3$& $h_4$& $h_5$&$\dim_{\R}
\SSS(\g)$\\ 
\hline
2&2&1&(0,0,12,13,14+23,34+52)& &\no  &\no &\no& \no\\
\extrarowheight10pt
2&2&1&(0,0,12,13,14,34+52)&  &\no &\no&\no&\no\\ 
2&3&1&(0,0,12,13,14,15)&&3&3&0& 7\\ 
2&3&1&(0,0,12,13,14+23,24+15)& &3,4 &2&0& 7\\
2&3&1&(0,0,12,13,14,23+15)& &2 &2&0& 7\\ 
2&4&2&(0,0,12,13,23,14)& &4 &4&0&8\\ 
2&4&2&(0,0,12,13,23,14-25)& &4 &2,3,4&0& 8\\
2&4&2&(0,0,12,13,23,14+25)& &4 &4&0& 8\\
\hline
3&4&2&(0,0,0,12,14-23,15+34)& &2 &2&0& 7\\ 
3&5&2&(0,0,0,12,14,15+23)& &4&4&2& 8\\ 
3&5&2&(0,0,0,12,14,15+23+24)& &4,5 &3,4&0,2& 8\\
3&5&2&(0,0,0,12,14,15+24)&1+5&5&4 &2 & 8\\ 
3&5&2&(0,0,0,12,14,15)&1+5&5&4&2 & 8\\ 
3&5&3&(0,0,0,12,13,14+35)& &\no &\no&\no&\no\\
3&5&3&(0,0,0,12,23,14+35)& &\no &\no&\no& \no\\ 
3&5&3&(0,0,0,12,23,14-35)& &\no&\no &\no &\no\\ 
3&5&3&(0,0,0,12,14,24)&1+5&\no&\no &\no &\no\\
3&5&3&(0,0,0,12,13+42,14+23)& &5 &3 &0 & 8\\ 
3&5&3&(0,0,0,12,14,13+42)& &5 &3 &0 & 8\\ 
3&5&3&(0,0,0,12,13+14,24)& &5 &2,3 &0 &8\\
3&6&3&(0,0,0,12,13,14+23)& &5,6,7 &3,4 &0 & 9\\ 
3&6&3&(0,0,0,12,13,24)& &5,6 &5 &0 & 9\\ 
3&6&3&(0,0,0,12,13,14)& &5,6 &4&0 & 9\\ 
3&8&4&(0,0,0,12,13,23)& &9,10&7,8&0 & 9\\
\hline 
4&6&3&(0,0,0,0,12,15+34)& &\no &\no&\no&\no\\
4&7&3&(0,0,0,0,12,15)&1\!+\!1\!+\!4&6&3&2& 9\\
4&7&3&(0,0,0,0,12,14+25)&1+5&6,7&3&2& 9\\ 
4&8&4&(0,0,0,0,13+42,14+23)& & 8&7&2& 10\\ 
4&8&4&(0,0,0,0,12,14+23)& & 8&6 &2 & 10\\
4&8&4&(0,0,0,0,12,34)&3+3& 8&7 &2  & 10\\ 
4&9&4&(0,0,0,0,12,13)&1+5& 10&7,8 &2 & 11\\
\hline 
5&9&4&(0,0,0,0,0,12+34)&1+5&\no&\no &\no & \no\\
5&11&4&(0,0,0,0,0,12)&1\!+\!1\!+\!1\!+\!3&13 & 9 &4 & 12\\
\hline
6&15&5&(0,0,0,0,0,0)&$1+\cdots+1$& 20 & 15 &6& 15\\
\hline 
\end{tabular}}
\end{center}\vfil\eject

\m Now we explain how to compute $h_3$. First, we have the following Lemma.

\begin{lemma}
\label{ker}
Let $(M,\omega)$ be a
symplectic manifold of dimension 6. Then
$$
h_3(M)=h_5(M)+\dim \ker
\{L\colon H^3(M)\longrightarrow H^5(M)\}.
$$
\end{lemma}

\p
Because of Theorem \ref{harm} and Proposition \ref{small},

\begin{equation}\label{uno}
H^3_{hr}(M)=P^3(M)+L(H^1(M)),
\end{equation}

where 

\begin{eqnarray*}
P^3(M)&=&\{v\in H^3(M)\mid v\wedge [\omega]=0\}\\ 
&=& \ker\{L\colon H^3(M)\longrightarrow H^5(M)\}.
\end{eqnarray*}

We need to compute the intersection $P^3(M)\cap L(H^1(M))$.
We set $A=\ker \{L^2\colon H^1(M)\longrightarrow H^5(M) \}$. Then

\begin{eqnarray*}
P^3(M)\cap L(H^1(M)) &=&\{a\wedge[\omega] \mid a\in H^1(M) \mbox{ and }
L^2(a)=0 \} \\
&=& {\rm Im}\, \{L_{\mid_A} \colon A\longrightarrow
H^3(M)\}. 
\end{eqnarray*}

Clearly, $\dim A=\dim \ker L_{\mid_A} + \dim(P^3(M)\cap
L(H^1(M))$. But 

\begin{eqnarray*}
\ker L_{\mid_A}&=&\{a\in H^1(M) \mid a\wedge[\omega]\\
&=&0 \mbox{ and } a\wedge[\omega]^2=0\} =\ker
\{L\colon H^1(M)\longrightarrow H^3(M)\}. 
\end{eqnarray*}

Thus
$$
\dim(P^3(M)\cap L(H^1(M)))= \dim A - \dim \ker \{L\colon
H^1(M)\longrightarrow H^3(M)\}.
$$

Taking into account that 
$$
\dim A + \dim L^2(H^1(M)) =
\dim \ker \{L\colon H^1\longrightarrow H^3\} + \dim L(H^1)=b_1,
$$
we conclude that 
\begin{equation}\label{dos}
\dim(P^3\cap L(H^1(M)))= \dim L(H^1(M)) - \dim L^2(H^1(M)).
\end{equation}

Since $h_5=\dim L^2(H^1(M))$, we deduce from (\ref{uno}),
(\ref{dos}) that 
$$ 
h_3=h_5+\dim \ker \{L\colon H^3(M)\longrightarrow
H^5(M)\}.
$$
\qed

Now, using the Nomizu Theorem, one can compute 
$$
\dim \ker \{L\colon H^3(M)\longrightarrow H^5(M)\}
$$
and therefore to compute $h_3$. As an example, we compute $h_3(M)$ where $M$ is 
the nilmanifold (0,0,0,12,14,15+23+24). Below we write $\ga_{ij\ldots k}$ 
instead of $\ga_i\wedge \ga_j \wedge \cdots \wedge\ga_k$.

First, by Nomizu theorem the cohomology groups of degrees 3 and 5
are: 
$$
\begin{array}{ccl}
H^3(M)&= \R^6 = &\span\{
[\alpha_{126}],[\alpha_{135}],
[\alpha_{136}+\alpha_{146}],[\alpha_{136}+\alpha_{235}],\\
\!\!\!\!&& [\alpha_{156}-\alpha_{236}-\alpha_{246}],
[\alpha_{156}+\alpha_{345}-\alpha_{246}]\},\\[5pt]
H^5(M)&= \R^3 = &\span\{ [\alpha_{12456}],[\alpha_{13456}],
[\alpha_{23456}]\}. 
\end{array}
$$
>From \cite{IRTU} we know that the cohomology class of any symplectic form
$\omega$ on $M$ must be a linear combination $$[\omega]=A\,
[\alpha_{13}]+B\, [\alpha_{15}]+C\, [\alpha_{23}]+ D\,
[\alpha_{16}+\alpha_{25}-\alpha_{34}]+ E\,
[\alpha_{26}-\alpha_{45}],$$ where $A,B,C,D,E \in \R$ and satisfy
$AE^2+BDE-CDE-D^3\not= 0$.

A direct calculation shows that the linear mapping $L\colon
H^3(M)\longrightarrow H^5(M)$ has, with respect to the bases of
$H^3(M)$ and $H^5(M)$ given above, the following matrix:

\[
\left(
\begin{array}{ccc}
-E & 0 & 0\\ 
0 & 0 & 0\\ 
-D & -E & 0\\
0 & -E & 0\\ 
-B & -D & E\\ 
-B & -2D & -E
\end{array}
\right)
\]

Therefore, if $E\not=0$ then $\dim \ker\{L\colon
H^3(M)\longrightarrow H^5(M)\}=3$; and if $E=0$ then this
dimension is 4, because the symplecticity condition implies that
$D\not=0$.

Finally, since we know from \cite{IRTU} that $h_4=4, h_5=2$ if $E\not=0$,
and $h_4=3, h_5=0$ if $E=0$, by Lemma \ref{ker} the following holds:
\begin{enumerate}
\item[{\rm (i)}] if $E\not=0$ then $h_3=5$;
\item[{\rm (ii)}] if $E=0$ then $h_3=4$.
\end{enumerate}

\m 
Because of Corollary \ref{flex3}, the nilmanifolds 
(0,0,12,13,23,14-25),\hfill\linebreak
(0,0,0,12,14,15+23+24), (0,0,0,12,13+14,24), 
(0,0,0,12,13,14+23),\hfill\linebreak 
(0,0,0,12,13,23), (0,0,0,0,12,13) are flexible. Now we prove that the manifolds
(0,0,12,13,14+23,24+15), (0,0,0,12,13,24),(0,0,0,12,13,14)\hfill\linebreak
and (0,0,0,0,12,14+25) are also flexible. So, summarizing, 
we have ten 6-dimensional flexible nilmanifolds. 

\begin{corollary}
\label{3-flex}
 Let $M$ be a closed $6$-dimensional manifold, and let $\omega, \omega'$ be two 
symplectic forms on $M$ such that $h_5(\omega)=h_5(\omega')$ and 
$h_3(\omega)>h_3(\omega')$. Then $M$ is $3$-flexible.
\end{corollary}

\p Let $k(\eta)=\dim\ker \{L_{\eta}: H^3(M) \to H^5(M)\}$ for a symplectic form 
$\eta$ on $M$. Asserting as in Corollary \ref{form}, we can prove that
$$
\gL:=\{\gl \in \R\bigm|k(\omega+\gl \omega')\le k(\omega)\}
$$
is an open and dense subset of $\R$. So, for every $\eps>0$, there exists 
$\gl\in (0,\eps)$ such that
$k(\omega+\gl\omega')>k(\omega)$. Now, choose $\gl$ so small that 
$\omega_t:=\omega+t\omega'$ is a 
symplectic form for $t\in [0,\gl]$. Because of what we said above, we have 
$k(\omega_0)=k(\omega)$ and $k(\omega)>k(\omega_{\gl})$. So, because of Lemma 
\ref{ker}, 
$$
h_3(\omega_0)>h_3(\omega_{\gl}),
$$   
and thus $M$ is 3-flexible.
\qed

Now we show an analog of Corollary \ref{strat} for $h_3$. It is interesting to 
notice that in this case the inequality in (iii) changes the direction, i.e. the 
generic numbers $h_3$ are minimal, unlike the case of Corollary \ref{strat}.

\begin{corollary}
\label{strat2}
Let $M^{6}$ be a manifold that admits a symplectic structure. Suppose that 
$h_{5}(M,\omega)$
does not depend on the symplectic structure $\omega$ on $M$. 
Then the following
three conditions are equivalent:
\par {\rm (i)} the set $\Omega (b, 3)$ is open and dense in $\osi(M)$;
\par {\rm (ii)} the interior of the set $\Omega (b, 3)$ in $\osi(M)$ is 
non-empty;
\par {\rm (iii)} the set $\Omega (b, 3)$ is non-empty and
$h_3(M,\omega)\geq b$ for every $\omega\in \osi(M)$.
\end{corollary}

\p This can be proved similarly to Corollary \ref{strat}, using Lemma \ref{ker}.
We explain why the inequality in (iii) changes the direction. It happens because 
we must use the 
fact that small perturbation does not {\it increase} the dimension of the kernel 
of a linear map (while it does not {\it decrease} the rank).  
\qed

\begin{comment}
\rm Here we show that the nilmanifold $M$ of the type
(0,0,0,12,14,15+23+24) is 3-flexible. Notice that this does not follow from  
Corollary \ref{3-flex}. Define
$$
\omega_t=(1-\cos t) \alpha_{13} -\cos t
(\alpha_{16}+\alpha_{25}-\alpha_{34}) + (1-\cos
t)(\alpha_{26}-\alpha_{45}).
$$ 
Each closed 2-form $\omega_t$ is
non-degenerate because the symplecticity condition
$AE^2+BDE-CDE-D^3= 3 \cos^2 t-3\cos t+1\not= 0$ for all $t\in \R$.
Thus, we have a complete family $\omega_t$ which is also
periodic, i.e. $\omega_{t+2\pi} =\omega_t$. Since $E=1-\cos t$,
this closed curve $\omega_t$ has the following properties:
\begin{enumerate}
\item[{\rm (i)}] $h_3(\omega_{2\pi k})=4$, $h_4(\omega_{2\pi k})=3$,
$h_5(\omega_{2\pi k})=0$, for any integer $k$;
\item[{\rm (ii)}] $h_3(\omega_t)=5$, $h_4(\omega_t)=4$,
$h_5(\omega_t)=2$, for $t\not= 2\pi k$.
\end{enumerate}
\end{comment}

\begin{remark}\rm
Sakane and Yamada \cite{SY} rediscovered some of our examples of flexible 
6-dimensional nilmanifolds.
\end{remark}

\section*{Acknowledgments} The first and fourth authors were partially
supported by the project UPV 127.310-EA7781/2000 and DGICYT
PB97-0504-C02-01/02. The second author was partially supported by the Fields 
Institute for Research in Mathematics Sciences, Toronto, Canada. The third 
author was partially financed by the Polish Research Committee (KBN).


\begin{thebibliography}{999}

\bibitem{AK} V.I. Arnold, B.A. Khesin: {\sl Topological Methods in 
Hydrodynamics},
Applied Math. Sci. 125, Springer, Berlin, 1998.

\bibitem{Br} J.-L. Brylinski: A differential complex for Poisson
manifolds,
{\sl J. Diff. Geom.} {\bf 28} (1988), 93-114.

\bibitem{GH} P. Griffiths, J. Harris: {\sl Principles of algebraic geometry}, 
Pure and Applied Mathematics. Wiley-Interscience [John Wiley \& Sons], New 
York, 
1978. 

\bibitem{IRTU} R. Ib\'a\~nez, Yu. Rudyak, A. Tralle, L. Ugarte: On 
symplectically harmonic forms on six-dimensional nilmanifolds. {\sl Comment. 
Math. Helv.} {\bf 76} (2001), no. 1, 89--109, Erratum 76 (2001), no 3, 576.

\bibitem{K} J.-L. Koszul: Crochet de Schouten-Nijenhuis et cohomologie,
The mathematical heritage of \'Elie Cartan (Lyon, 1984). 
{\sl AstŽrisque} 1985, Numero Hors Serie, 257--271. 

\bibitem{L} P. Liberman: Sur le probl\`eme d'\'equivalence de certaines 
structures infinit\'esimales r\'egulieres, {\sl Ann. Mat. Pura Appl.} {\bf 36} 
(1954), 27--120. 

\bibitem{LM} P. Libermann, C. Marle: {\sl Symplectic Geometry and Analytical 
Mechanics},
Kluwer, Dordrecht, 1987.

\bibitem{M} A.I. Malcev: On a class of homogeneous spaces, {\sl Izv AN SSSR 
Ser. 
Matem.} {\bf 3} (1949),9-32.

\bibitem{Ma} O. Mathieu: Harmonic cohomology classes of symplectic
manifolds,
{\sl Comment. Math. Helvetici} {\bf 70} (1995), 1-9.

\bibitem{OV} A.L. Onishchik, E.B. Vinberg (Eds): Lie Groups and Lie Algebras,
Encycl. Math. Sci., 41, Springer, Berlin, 1994.

\bibitem{R} M. Raghunathan: {\sl Discrete subgroups of Lie groups}, Springer, 
Berlin, 
1972

\bibitem{SY} Y. Sakane, T. Yamada: Harmonic cohomology groups of compact 
symplectic nilmanifolds, to appear in Proc. of 9-th MSJ-IRI.

\bibitem{Sa} S. Salamon: Complex structures on nilpotent Lie algebras, {\sl J. 
of Pure and Appl. Algebra} {\bf 157} (2001), no. 2-3, 311--333.

\bibitem{TO} A. Tralle, J. Oprea: {\sl Symplectic Manifolds with no
K\"ahler
Structure}, Lecture Notes in Math. {\bf 1661}, Springer, Berlin 1997.

\bibitem{W}
R. Wells, Differential analysis on complex manifolds. Second edition. Graduate 
Texts in Mathematics, 65. Springer-Verlag, New York-Berlin, 1980. 

\bibitem{Ym} 
T. Yamada: Harmonic cohomology groups of compact symplectic 
nilmanifolds, to appear in Osaka J. Math. 

\bibitem{Y} 
D. Yan: Hodge structure on symplectic manifolds, {\sl Adv.
in Math.} {\bf 120} (1996), 143-154.

\end{thebibliography}
\end{document}